\documentclass[11pt,a4paper]{article}
\usepackage{indentfirst}
\usepackage{amsmath,amssymb,amsfonts,amsthm,graphics}
\usepackage{makeidx}
\usepackage{color}

\DeclareMathOperator{\asc}{asc}
\DeclareMathOperator{\rlm}{rlm}
\newcommand{\cA}{\mathcal{A}}
\newcommand{\ree}[1]{(\ref{#1})}
\newcommand{\beq}{\begin{equation}}
\newcommand{\eeq}{\end{equation}}
\newcommand{\bea}{\begin{eqnarray*}}
\newcommand{\eea}{\end{eqnarray*}}
\newcommand{\si}{\sigma}
\newcommand{\ep}{\epsilon}

\newtheorem{theorem}{Theorem}[section]

\newtheorem{conjecture}{Conjecture}[section]
\numberwithin{equation}{section}
\begin{document}
\title{\Large\bf On 021-Avoiding Ascent Sequences}
\author{William Y.C. Chen\\[-5pt]
\small Center for Combinatorics, LPMC-TJKLC, Nankai University\\[-5pt]
\small Tianjin 300071, P.R. China, {\tt chen@nankai.edu.cn}
\\[5pt]
Alvin Y.L. Dai\\[-5pt]
\small Center for Combinatorics, LPMC-TJKLC, Nankai University\\[-5pt]
\small Tianjin 300071, P.R. China, {\tt alvin@cfc.nankai.edu.cn}
\\[5pt]
Theodore Dokos\\[-5pt]
\small Department of Mathematics, University of California, Los Angeles\\[-5pt]
\small Box 951555, Los Angeles, CA 90095-1555, USA, {\tt t.dokos@gmail.com}
\\[5pt]
Tim Dwyer\\[-5pt]
\small Department of Mathematics,  Dartmouth College,\\[-5pt]
\small Hanover, NH 03755, USA, {\tt jtimdwyer@gmail.com}
\\[5pt]
Bruce E. Sagan\\[-5pt]
\small Department of Mathematics, Michigan State University,\\[-5pt]
\small East Lansing, MI 48824-1027, USA, {\tt sagan@math.msu.edu}
}
\date{}
\maketitle

\begin{abstract}
Ascent sequences were introduced by Bousquet-M\'{e}lou,
Claesson,  Dukes and Kitaev in their study of $(\bf{2+2})$-free posets. An ascent sequence of length $n$ is a nonnegative integer sequence $x=x_{1}x_{2}\ldots x_{n}$ such that $x_{1}=0$ and
$x_{i}\leq \asc(x_{1}x_{2}\ldots x_{i-1})+1$ for all $1<i\leq n$, where $\asc(x_{1}x_{2}\ldots x_{i-1})$ is the number of ascents in the sequence $x_{1}x_{2}\ldots x_{i-1}$.
We let $\cA_n$ stand for the set of such sequences and use $\cA_n(p)$ for the subset of sequences avoiding a pattern $p$.
Similarly, we let $S_{n}(\tau)$ be the set of $\tau$-avoiding permutations in the symmetric group $S_{n}$.  Duncan and Steingr\'{\i}msson have shown that the ascent statistic has the same distribution over  $\cA_n(021)$ as over $S_n(132)$.  Furthermore, they conjectured that the pair  $(\asc, \rlm)$ is equidistributed over $\cA_n(021)$ and  $S_n(132)$ where $\rlm$ is the right-to-left minima statistic.  We prove this conjecture by constructing a bistatistic-preserving bijection.

\vskip 2mm

\noindent
{\bf Keywords}: $021$-avoiding ascent sequence, $132$-avoiding permutation, right-to-left minimum, number of ascents,  bijection \\[2mm]

\noindent
{\bf AMS subject classification:} 05A05, 05A19
\end{abstract}

\section{Introduction}

 The objective of this note is to establish a bijection
 which leads to the equidistribution of the
 pair of statistics $(\asc, \rlm)$ over $021$-avoiding ascent sequences and over $132$-avoiding permutations. This confirms a conjecture posed by Duncan and Steingr\'{\i}msson \cite{P.D}.

Let us give an overview of the notation and terminology.
Let $S_{n}$ denote the set of permutations of $[n]$, where $[n]=\{ 1, 2, \ldots, n\}$. Given a permutation $\pi\in S_{n}$ and a permutation $\tau\in S_{k}$,  we say that a subsequence $\pi_{i_{1}}\pi_{i_{2}}\ldots\pi_{i_{k}}$, $1\leq i_{1}<i_{2}<\cdots<i_{k}\leq n$, of $\pi$ is
of \emph{pattern} $\tau$ if it is order isomorphic to $\tau$, that is, this subsequence has the same relative order as $\tau$.
If $\pi$ does not contain any subsequence of pattern $\tau$, then we say that $\pi$ \emph{avoids} $\tau$, or $\pi$ is \emph{$\tau$-avoiding}. We denote by $S_{n}(\tau)$
the set of $\tau$-avoiding permutations in $S_{n}$. For example, the permutation $763894512$ contains the subsequence $3952$ of pattern $2431$, but it is $1234$-avoiding. Pattern avoiding permutations have been intensively studied in recent years from many points of view, see \cite{M.B, C.K, R.S}.

Ascent sequences were introduced by Bousquet-M\'{e}lou, Claesson,  Dukes and Kitaev \cite{B.M}. For a sequence $x=x_{1}x_{2}\ldots x_{n}$ of nonnegative integers, we say that an index $i$ $(1\leq i<n)$ is an \emph{ascent} if $x_{i}<x_{i+1}$. We denote by $\asc(x)$ or merely $\asc x$ the number of ascents of $x$.  A sequence $x=x_{1}x_{2}\ldots x_{n}$ is called an \emph{ascent sequence} if $x_{1}=0$ and
$$
x_{i}\leq \asc(x_{1}x_{2}\ldots x_{i-1})+1
$$
for all $1<i\leq n$.   For example, $x=010122$ is an ascent sequence while $x=010142$ is not since $x_5=4>\asc(0101)+1=3$.
We let $\cA_n$ denote the set of ascent sequences of length $n$.
For an ascent sequence, a \emph{pattern} is a word on a nonnegative integers $\{0, 1, \ldots, k\}$, where each element $i$ appears at least once.  Containment and avoidance of patterns for ascent sequences are defined in the same way as for permutations.
For example, the ascent sequence $01231234$ has five occurrences of the pattern $001$, namely, the subsequences $112, 113, 114, 223, 224$,
and the ascent sequence $01012203$ is $021$-avoiding. We denote by $\cA_n(p)$ the set of ascent sequences of length $n$ avoiding pattern $p$.

In addition to the ascent statistic, we will be interested in the number of right-to-left minima.  A \emph{right-to-left minimum} of any seqence $x$ of nonnegative integers is  an index $i$ such that $x_{i}<x_{j}$ for all $j>i$. The number of right-to-left minima of $x$ is denoted by $\rlm(x)=\rlm x$. For example, $\rlm(010122)=3$.

Ascent sequences are closely connected to   $(\bf{2+2})$-free posets~\cite{B.M}, upper-triangular matrices~\cite{M.D}, Stoimenow's matchings~\cite{A.C},
and the Catalan numbers $C_n$~\cite{P.D}.
In particular£¬a poset is called  \emph{$(\bf{2+2})$-free} if it does not contain an induced subposet which is isomprphic to the disjoint union of two $2$-element chains. Bousquet-M\'{e}lou, Claesson, Dukes and Kitaev~\cite{B.M} found a bijection from $(\bf{2+2})$-free posets to ascent sequences which maps the number of levels of the poset to the number of ascents of the sequence.
Dukes and Parviainen~\cite{M.D} established a bijection between ascent sequences and nonnegative upper-triangular matrices.
Duncan and  Steingr\'{\i}msson~\cite{P.D}
have shown that $\#\cA_n(p)=C_n$ for any of the patterns $p=101$, $0101$, or $021$. It is well known that $\#S_n(132)=C_n$ and Duncan-Steingr\'{\i}msson also showed that the ascent statistic is equidistributed over $\cA_n(021)$ and $S_n(132)$. Furthermore, they proposed the following conjecture.

\begin{conjecture}\label{con1}
The bistatistic $(\asc, \rlm)$ has the same distribution over $\cA_n(021)$ and $S_n(132)$.
\end{conjecture}

The objective of this note is to give a bijective proof of  the above conjecture.

\section{Proof of the conjecture}

In order to construct our bijection, we will need the concept of the special maximum value of an ascent
sequence $x=x_1 x_2\ldots x_n$.  The \emph{special maxium value} of $x$ is the largest integer $M$ such that there is an index $i$ with $x_i=M$ and
\beq
\label{special}
x_{i}= \asc(x_{1}x_{2}\ldots x_{i-1})+1
\eeq
giving us equality in the defining relation for an ascent sequence. To illustrate the notion, if $x=01013312434$ then $M=3$.  Note that, except for the zero sequence, $M$ will always exist since the first $1$ in any nonzero sequence satisfies~\ree{special}.  So we define $M=0$ for a zero sequence.

Also define a \emph{special maximum index} as an index $i$ where $x_i$ satisfies $x_i=M$ as well as condition~\ree{special}.  The first index $i$ with $x_i=M$ is always a special maximum index.  In fact, if $[i,j]$ is the largest interval of indices starting with the first special maximum index and satisfying $x_i=x_{i+1}=\dots=x_j=M$ then we claim that these are exactly the special maximum indices.  To see this, first note that if $k\in[i,j]$ then $k$ is a special maximum index because
$$
x_k=x_i= \asc(x_{1}x_{2}\ldots x_{i-1})+1= \asc(x_{1}x_{2}\ldots x_{k-1})+1
$$
since there are no ascents between $x_i$ and $x_k$.  To see that no other index can be special maximum, suppose $x_k=M$ with $k\ge j+2$.  Now $x_j>x_{j+1}$ because if $x_j<x_{j+1}$ then the special maximum value would be at least $M+1$.  Thus there must be an ascent between $x_{j+1}$ and $x_k$ so that~\ree{special} is no longer an equality when $i=k$.  Call the special maximum value \emph{unique} if there is only one special maximum index and \emph{repeated} otherwise.

\begin{theorem}
The bistatistic $(\asc, \rlm)$ has the same distribution over $\cA_n(021)$ and $S_n(132)$.
\end{theorem}
\noindent{\bf Proof.}
We will inductively build a bijection  $\phi_n\colon\cA_n(021)\rightarrow S_n(132)$ preserving the bistatistic.  To do so, we will need decompositions of $\cA_n(021)$ and $S_n(132)$ into pieces indexed by smaller subscripts.  We will start on the ascent side.

A simple but  important observation for what follows is  that $p\in\cA_n(021)$ if and only if the nonzero entries of $p$ are weakly increasing.  We will use this fact to construct a bijection $f=f_n$ between $\cA_n(021)$ and the set of pairs
$$
\bigcup_{i=1}^n \cA_{i-1}(021)\times \cA_{n-i}(021).
$$

Consider  $x\in\cA_n(021)$ and suppose first that $x$ has a repeated special maximum value $M$.  Let $k$ be any of the special maximum indices and define
$$
f(x)=(\epsilon,z)
$$
where $\epsilon$ is the empty sequence and $z$ is $x$ with $x_k$ removed.  For example, if $x=01013300304$ then $z=0101300304$.  Clearly $z$ still avoides $021$ since its nonzero entries still increase and, since $x$'s special maximum value was repeated, $z$ still has $M$ as its special maximum value.  Since the special maximum value does not change, one can construct an inverse map from $\cA_0(021)\times\cA_{n-1}(021)$ back to the elements of $\cA_n$ which have a repeated special maximum in the obvious way.  Finally note that in this case
$$
\asc x = \asc z  \quad\text{and}\quad  \rlm x = \rlm z.
$$

Now suppose that $x$ has a unique special maximum value $x_i=M$.  Here we let
$$
f(x)=(y,z)
$$
where $y=x_1\ldots x_{i-1}$ and $z$ is obtained from $x'=x_{i+1}\ldots x_n$ by subtracting $M-1$ from all the nonzero entries.  To illustrate, if $x=0101300304$ then $y=0101$ and $z=00102$.  It is clear that $y\in\cA_{i-1}(021)$.  To show that $z\in\cA_{n-i}(021)$, we first note that $z$ still has weakly increasing nonzero elements and so avoids $012$.  We must also demonstrate that $z$ is an ascent sequence.  Since the defining condition for an ascent sequence is trivial for zero elements, we need only consider $z_r\neq0$.  But since we have subtracted the same amount from all nonzero entries of $x'$,  the index $r$ is an ascent of $z$ if and only if the index $r+i$ is an ascent of $x'$.  Also, since $M$ was the special maximum
value, we have $\asc(x_1\ldots x_i) = M$, $x_i>x_{i+1}$ and $x_{r+i}\le\asc(x_1x_2\ldots x_{r+i-1})$ for any $r\ge1$.  It follows that for any $z_r\neq 0$ we have
\bea
z_r &=& x_{r+i}-M+1\\
&\le&\asc(x_1\ldots x_{r+i-1})-M+1\\
&=&\asc(x_1\ldots x_i)+\asc(x_{i+1}\ldots x_{r+i-1})-M+1\\
&=&\asc(z_1\ldots z_{r-1})+1
\eea
which is what we wished to prove.  Constructing the inverse of this part of the map is similar to what was done in the first case.

It will be useful to record what happens to our two statistics in the second case defining $f$.  For the ascent statistic, we have everything in place from the previous paragraph and the fact that, by definition of $M$, $x_i>x_{i-1}$.  Thus
\bea
\asc x&=&\asc(x_1\ldots x_i)+\asc(x_{i+1}\ldots x_n)\\
&=&\asc(x_1\ldots x_{i-1})+1+\asc(x_{i+1}\ldots x_n)\\
&=&\asc y + \asc z + 1.
\eea
In terms of right-to-left minima, we distinguish two subcases.  If $i\le n-1$ then, since $x_{i+1}=0$, the right-to-left minima
of $x$ must occur in the sequence $x_{i+1}\ldots x_n$.  Since the subtraction of $M-1$ does not change the positions of these minima, we have
$$
\rlm x = \rlm z.
$$
On the other hand, if $i=n$ then $z=\epsilon$ and $x_n$ is a right-to-left minimum, giving
$$
\rlm x= \rlm y +1.
$$

We will now review the standard decomposition of $S_n(132)$ which gives a bijection $g$ from this  set to
$$
\bigcup_{i=1}^n S_{i-1}(132)\times S_{n-i}(132).
$$
If $\pi\in S_n(132)$ then we write $\pi=\pi_L n\pi_R$ where $\pi_L,\pi_R$ are the elements to the left and right of $n$, respectively.  Define the index $i$ by $\pi_i=n$.  Then it is well known that $\pi\in S_n(132)$ if and only if $\pi_L, \pi_R$ avoid $132$ and every element of $\pi_L$ is bigger than every element of $\pi_R$.  So we let
$$
g(\pi)=(\rho,\si)
$$
where $\rho\in S_{i-1}(132)$ and $\si\in S_{n-i}(132)$ are order isomorphic to $\pi_L$ and $\pi_R$, respectively.

As with the ascent sequence decomposition, we have to consider what happens to our statistics in two separate cases.  The first is when $\rho=\epsilon$, equivalently, $i=1$.  So $\pi=n\pi_R$ and so $n$ makes no contribution either to the ascents or right-to-left maxima.  It follows that
$$
\asc \pi=\asc\si  \quad\text{and}\quad  \rlm \pi = \rlm \si
$$
just as in the corresponding case for ascent sequences.

Now suppose $1<i\le n$.  So there will be an ascent ending at $n$ and all other ascents of $\pi$ correspond to ascents of $\rho$ or ascents of $\si$.  It follows that
$$
\asc\pi=\asc\rho+\asc\si+1.
$$
For the right-to-left minima we again break into two subcases depending on whether the second component of our bijection is $\epsilon$ or not.  If $\si\neq\epsilon$ then $i<n$ and the right-to-left minima of $\pi$ are all in $\pi_R$ because of the relative sizes of the elements of $\pi_R$ and $\pi_L$.  This gives
$$
\rlm \pi=\rlm \si.
$$
Now consider $\si=\ep$ so that $\pi=\pi_Ln$.  Thus $n$ is a right-to-left minimum of $\pi$ as is every right-to-left minimum of $\pi_L$.  So in this subcase
$$
\rlm \pi=\rlm\rho + 1.
$$

Finally, we construct  $\phi_n\colon\cA_n(021)\rightarrow S_n(132)$ as follows.  Start with $\phi_0(\epsilon)=\epsilon$.  Assuming that $\phi_i$ has been defined for all $i<n$, we define $\phi_n$ to be the composition
$$
\cA_n(021)
\stackrel{f}{\rightarrow}
\bigcup_{i=1}^n \cA_{i-1}(021)\times \cA_{n-i}(021)
\stackrel{h}{\rightarrow}
\bigcup_{i=1}^n S_{i-1}(132)\times S_{n-i}(132)
\stackrel{g^{-1}}{\rightarrow}
S_n(132)
$$
where the restriction of $h$ to  $\cA_{i-1}(021)\times \cA_{n-i}(021)$ is $\phi_{i-1}\times\phi_{n-i}$.  It should be clear from the equations derived for $\asc$ and $\rlm$ when defining $f$ and $g$ that this bijection preserves the bistatistic.
\hfill\qed

\vskip 3mm

\noindent {\bf Acknowledgments.}
We wish to thank Einar Steingr\'{\i}msson for helpful comments.
This work was supported by the 973 Project, the PCSIRT Project
of the Ministry of Education, and the National Science Foundation of China.


\begin{thebibliography}{20}

\bibitem{M.B}
M. B\'{o}na, Combinatorics of Permutations, CRC Press, 2004.

\bibitem{B.M}
M. Bousquet-M\'{e}lou, A. Claesson, M. Dukes and S. Kitaev,
(2+2)-free posets, ascent sequence and pattern avoiding
permutations, J. Combin. Theory, Ser. A, 117 (7) (2010), 884-909.

\bibitem{A.C}
A. Claesson, M. Dukes and S. Kitaev, A direct encoding of
 Stoimenow's matchings as ascent sequences, Australas.
 J. Combin., 49 (2011), 47-59.

\bibitem{M.D}
M. Dukes and R. Parviainen,  Ascent sequences and upper
 triangular matrices containing non-negative integers,
 Electron. J. Combin., 17 (1) (2010), \#R53.

\bibitem{P.D}
P. Duncan and E. Steingr\'{\i}msson, Pattern avoidance
in ascent sequences, Electron. J. Combin.,  18 (2011), \#P226.

\bibitem{C.K}
C. Krattenthaler, Permutations with restricted patterns
and Dyck paths, Adv. Appl. Math., 27 (2001), 510-530.

\bibitem{R.S}
R. Simion and F.W. Schmidt, Restricted permutations,
European J. Combin., 6 (1985), 383-406.

\end{thebibliography}
\end{document}